\documentclass[12pt]{article}

\input{tcilatex}
\begin{document}

\title{A general decay and optimal decay result in a heat system with a
viscoelastic term}
\author{Abderrahmane Youkana$^{(1)}$, Salim A. Messaoudi$^{(2)}$ \& Aissa Guesmia$%
^{(3)}$ \\
(1) Department of Mathematics, University of Batna 2 \\
Batna 05078, Algeria\\
E-mail: abderrahmane.youkana@univ-batna2.dz\\
\qquad  \ \ abder.youkana@yahoo.fr\\
(2) Department of Mathematics and Statistics\\
KFUPM, Dhahran 31261\\
Saudi Arabia\\
E-mail: messaoud@kfupm.edu.sa\\
(3) Elie Cartan Institute of Lorraine, Bat. A \\
Lorraine - Metz University, Ile de Saulcy, 57045\\
Metz Cedex 01, France\\
E-mail: aissa.guesmia@univ-lorraine.fr}
\date{}
\maketitle

\begin{abstract}
We consider a quasilinear heat system in the presence of an integral term
and establish a general and optimal decay result from which improves and
generalizes several stability results in the literature.
\end{abstract}

\section{Introduction}

In this work, we consider the following problem: 
\begin{equation}
\left\{ 
\begin{array}{ll}
A(t)|u_{t}|^{m-2}u_{t}-\Delta u+\int_{0}^{t}g(t-s)\Delta u(x,s)\ ds=0 & 
\text{in $\Omega \times (0,+\infty ),$} \\ 
u(x,t)=0 & \text{in }\partial \text{$\Omega \times \mathrm{I\hskip-2ptR}^{+},$} \\ 
u(x,0)=u_{0}(x) & \text{in $\Omega ,$}
\end{array}
\right. 
\end{equation}
where $m\geq 2$,\ $\Omega $ is a bounded domain of $\mathrm{I\hskip-2ptR}^{n}
$, $n\in \mathrm{I\hskip-2ptN}^{*}$, with a smooth boundary $\partial \Omega $,\ $g:\,\mathrm{I\hskip-2ptR}^{+} \rightarrow \mathrm{I\hskip-2ptR}^{+}$ is a positive nonincreasing function, and $A:\,\mathrm{I\hskip-2ptR}^{+} \rightarrow M_n (\mathrm{I\hskip-2ptR})$ is a bounded square matrix satisfying $A\in C(\mathrm{I\hskip-2ptR}^{+})$ and, for some positive constant $c_{0},$ 
\begin{equation}
\left( A(t)v,v\right) \geq c_{0}|v|^{2},\qquad \forall t\in \mathrm{I\hskip%
-2ptR}^{+},\ \forall v\in \mathrm{I\hskip-2ptR}^{n},
\end{equation}
where\ $(.,.)$ and $|.|$ are the inner product and the norm, respectively,
in $\mathrm{I\hskip-2ptR}^{n}$. The equation in consideration arises from
various mathematical models in engineering and physics. For instance, in the
study of heat conduction in materials with memory, the classical Fourier law
is replaced by the following form (cf. $\left[ 9\right] )$: 
\[
\mbox{ }q=-d\nabla u-\int_{-\infty }^{t}\nabla \left( k(x,t)u(x,\tau
)\right) \ d\tau ,
\]
where $u$ is the temperature, $d$ the diffusion coefficient and the integral
term represents the memory effect in the material. This type of problems has
considered by a number of researchers; see $\left[ 2,9,11\right] $ and the
references therein$.$ From a mathematical point of view, we expect that the
integral term would be dominated by the leading term in the equation, so
that the theory of parabolic equation can be applied. In fact, this has been
confirmed by the work of Yin $\left[ 11\right] ,$ in which he considered a
general equation of the form 
\[
u_{t}=divA(x,t,u,u_{x})+a(x,t,u,u_{x})+\int_{0}^{t}divB(x,t,\tau
,u,u_{x})\ d\tau 
\]
and proved the existence of a unique weak solution under suitable conditions
on $A,$ $B$ and $a$. See more results concerning global existence and
asymptotic behavior in Nakao and Ohara $[7]$, Nakao and Chen $[8]$, and
Engler \textit{et al}. $[3]$. Pucci and Serrin $[10]$ discussed the
following system: 
\[
A(t)|u_{t}|^{m-2}u_{t}=\Delta u-f(x,u),
\]
for $m>1$ and $f$ satisfying 
\[
(f(x,u),u)\geq 0
\]
and showed that strong solutions tend to the rest state as $t\rightarrow +\infty ,$ however, no
rate of decay has been given. Berrimi and Messaoudi $\left[ 1\right] $
showed that, if $A$ satisfies (1.2), then solutions with small initial
energy decay exponentially for $m=2$ and polynomially if $m>2.$ Messaoudi
and Tellab [5] considered (1.1), under condition (1.2) and for relaxation
function $g$ satisfying a general decay condition of the form 
\[
g^{\prime}(t)\leq -\xi (t)g(t),\quad\forall t\in\mathrm{I\hskip-2ptR}^{+},
\]
for some nonincreasing differentiable function $\xi :\mathrm{I\hskip-2ptR}^{+}\rightarrow \mathrm{I%
\hskip-2ptR}^{+}$, and established a general decay result, from which
the exponential and polynomial decay rates of [1] are only special cases.
Recently, Liu and Chen [4] investigated (1.1), with a nonlinear source term,
and established a general decay result under suitable conditions on $g$\ and
the nonlinear source term. They also proved a blow-up result for the
solution with both positive and negative initial energy.\newline

In this work, we discuss (1.1) when $g$ is of a more general decay, and
establish a general and optimal decay result, which improves those of
Berrimi and Messaoudi [1], Liu and Chen [4], and Messaoudi and Tellab [5].

\section{Preliminaries}

\setcounter{equation}{0}In this section, we present some material needed in
the proof of our result. For the relaxation function $g$ we assume that\vspace{%
0.5cm}\newline
$(G_{1})$ The function $g:\,\mathrm{I\hskip-2ptR}^{+}\rightarrow \mathrm{I%
\hskip-2ptR}^{+}$ is a differentiable function satisfying 
\[
g(0)>0\quad\mbox{ }\mbox{and}\quad 1-\int_{0}^{+\infty }g(s)\ ds=l>0.
\]
$(G_{2})$ There exist a constant $p\in [1,3/2)$ and a nonincreasing
differentiable function $\xi :\mathrm{I\hskip-2ptR}^{+}\rightarrow \mathrm{I%
\hskip-2ptR}^{+}$ satisfying 
\[
g^{\prime }(t)\leq -\xi (t)g^{p}(t),\quad\forall t\in\mathrm{I\hskip-2ptR}^{+}.\newline
\]
$(G_{3})$ We also assume that 
\begin{eqnarray*}
2 &\leq &m\leq \frac{2n}{n-2} \quad\quad \mbox{if}\ n\geq 3, \\
m &\geq &2\quad \quad \quad \quad \quad\quad\,\,\mbox{if}\ n=1,2. \\
&&
\end{eqnarray*}
\textbf{Remark 2.1. }There are many functions satisfying (G1) and (G2).
Examples of such functions are, for $b>0,$ $\alpha >0$, $\nu >1,$ and $a>0$ small enough, 
\[
g_{1}(t)=ae^{-b(t+1)^{\alpha }}\quad \mbox{ }\mbox{and}\quad \mbox{ }g_{2}(t)=%
\frac{a}{(1+t)^{\nu }}.
\]
\qquad We will also be using the embedding $H_{0}^{1}(\Omega
)\hookrightarrow L^{q}(\Omega ),$ $L^{r}(\Omega )\hookrightarrow
L^{q}(\Omega )$, for $2\leq q\leq r<+\infty ,$ and Poincar\'{e}'s
inequality. The same embedding constant $C_{*}$ will be used, and $C$ denotes a generic positive constant.\vspace{0.5cm}

We introduce the following: 
\begin{equation}
E(t)=\frac{1}{2}(g\circ \nabla u)(t)+\frac{1}{2}\left(
1-\int_{0}^{t}g(s)\ ds\right) ||\nabla u(t)||_{2}^{2} ,\quad\forall t\in\mathrm{I\hskip-2ptR}^{+},
\end{equation}
where $||.||_{q}=$ $||.||_{(L^{q}(\Omega ))^{n}},$ for $1\leq q<+\infty$, and 
\begin{equation}
(g\circ \nabla u)(t)=\int_{0}^{t}g(t-\tau )||\nabla u(.,t)-\nabla u(.,\tau
)||_{2}^{2}\ d\tau ,\quad\forall t\in\mathrm{I\hskip-2ptR}^{+}.
\end{equation}
Similarly to [10], we give the definition of a strong solution of (1.1).%
\vspace{0.5cm}\newline
\textbf{Definition 2.1.} A weak solution of (1.1) on $[0,T]$ is a function 
\[
u\in C\left( [0,T);(H_{0}^{1}(\Omega ))^{n}\right) \cap C^{1}\left(
(0,T);(L^{m}(\Omega ))^{n}\right) 
\]
which satisfies 
\[
\int_{0}^{t}\int_{\Omega }\left( \nabla u(x,s)-\int_{0}^{s}\nabla u(x,\tau
)d\tau \right) \cdot \nabla \phi (x,s)\ dx\ ds
\]
\[
+\int_{0}^{t}A(s)|u_{t}|^{m-2}u_{t}(x,s)\cdot \phi (x,s)\ dx\ ds=0,
\]
for all $t$ in $[0,T)$ and all $\phi $ in $C\left( [0,T);(H_{0}^{1}(\Omega
))^{n}\right) $.\vspace{0.5cm}\newline
\textbf{Remark 2.2.} Similarly to [10], we assume the existence of a
solution. For the linear case ($m=2$)$,$ one can easily establish the
existence of a weak solution by the Galerkin method. In the one-dimensional
case ($n=1$), the existence is established in a more general setting by Yin $%
\left[ 11\right] .\vspace{0.5cm}$

Finally, we state an important lemma [6].\vspace{0.5cm}\newline
\textbf{Lemma 2.1}. \emph{Assume that }$g$\emph{\ satisfies (}$G1$\emph{)
and (}$G2$\emph{) and }$u$\emph{\ is the solution of }$(1.1),$\emph{\ then} 
\emph{there exists a positive constant }$k_0$\emph{\ such that} 
\begin{equation}
\xi (t)(g\circ \nabla u)(t)\leq k_0\left( -E^{\prime }(t)\right)^{\frac{1}{%
2p-1}},\quad\forall t\in\mathrm{I\hskip-2ptR}^{+}.
\end{equation}

We also recall the following particular case of the well-known Jensen inequality which will be of essential
use in obtaining our result: let $f:\Omega \rightarrow \mathrm{I\hskip-2ptR}%
^{+}$ and $h:\Omega \rightarrow \mathrm{I\hskip-2ptR}^{+}$ be integrable
functions on $\Omega$ such that 
\[
\int\limits_{\Omega }h(x)\ dx=k>0.
\]
Then, for any $p>1$, we have 
\begin{equation}
\frac{1}{k}\int\limits_{\Omega }(f(x))^{\frac{1}{p}}h(x)\ dx\leq \left( \frac{1%
}{k}\int\limits_{\Omega }f(x)h(x)\ dx\right) ^{\frac{1}{p}}.
\end{equation}

\section{Decay result}

\setcounter{equation}{0}In this section, we state and prove our main result.
We start with a lemma. \vspace{0.5cm}\newline
\textbf{Lemma 3.1}. \emph{Let }$u$\emph{\ be the solution of }$(1.1).$\emph{%
\ Then the energy satisfies} 
\begin{equation}
E^{\prime }(t)=-\int_{\Omega }A(t)|u_{t}|^{m}\ dx-\frac{1}{2}g(t)||\nabla
u(t)||_{2}^{2}+\frac{1}{2}(g^{^{\prime }}\circ \nabla u)(t)\leq 0,\quad\forall t\in\mathrm{I\hskip-2ptR}^{+}.
\end{equation}
\textbf{Proof}. By multiplying the first equation in $(1.1)$ by $u_{t}$,
integrating over $\Omega $ we get (3.1), after routine manipulations.\vspace{%
0.5cm}\newline
\textbf{Lemma 3.2}. \emph{Let }$u$\emph{\ be a solution of problem }$(1.1)$%
\emph{. Then, for any \ }$\delta >0,$\emph{\ we have } 
\begin{equation}
\Vert \nabla u(t)\Vert _{2}^{2}\leq c_{4}\delta E(t)-\frac{C_{\delta }}{c_{0}%
}E^{\prime }(t)+c_{5}(g\circ \nabla u)(t), \quad\forall t\in\mathrm{I\hskip-2ptR}^{+} ,
\end{equation}
where $c_{0}$ is introduced in (1.2), $c_{4}$ and $c_{5}$ are two
positive constants, and $C_{\delta }$ is a positive constant depending on $\delta $. 
\vspace{0.5cm}\newline
\textbf{Proof}. Multiplying the first equation in $(1.1)$ by $u$ and integrating
over $\Omega ,$ we get 
\begin{equation}
\Vert \nabla u(t)\Vert _{2}^{2}=-\int_{\Omega }A(t)|u_{t}|^{m-2}u_{t}u(x,t)\
dx+\int_{\Omega }\int_{0}^{t}g(t-s)\nabla u(x,s)\cdot \nabla u(x,t)\ ds\ dx.
\end{equation}
Now, we estimate the right-hand side of (3.3). By using Young's and Poincar%
\'{e}'s inequalities, the boundedness of $A$, conditions $(G_{1})$ and $(G_{3})$, and the fact that 
\[
E(t)\leq E(0),
\]
we find, for any $\delta >0$, 
\begin{eqnarray}
-\int_{\Omega }A(t)|u_{t}|^{m-2}u_{t}u\ dx &\leq &\delta \Vert u(.,t)\Vert
_{m}^{m}+C_{\delta }\Vert u_{t}(.,t)\Vert _{m}^{m}  \nonumber  \label{8} \\
&\leq &\delta C_{*}^{m}\Vert \nabla u(.,t)\Vert _{2}^{m}+C_{\delta }\Vert
u_{t}(.,t)\Vert _{m}^{m}  \nonumber \\
&\leq &\delta C_{*}^{m}\left( \frac{2E(0)}{l}\right) ^{\frac{m-2}{2}}\left( 
\frac{2}{l}E(t)\right) +C_{\delta }\Vert u_{t}(.,t)\Vert _{m}^{m}  \nonumber
\\
&\leq &c_{1}\delta E(t)-\frac{C_{\delta }}{c_{0}}E^{\prime }(t).
\end{eqnarray}
Next, we estimate the second term of the right-hand side of (3.3) carefully.
By Young's inequality, we easily see that 
\begin{eqnarray}
&&\int_{\Omega }\nabla u(x,t)\ .\int_{0}^{t}g(t-s)\nabla u(x,s)\ ds\ dx\leq 
\frac{1}{2}\Vert \nabla u(.,t)\Vert _{2}^{2}  \nonumber \\
&+&\frac{1}{2}\int_{\Omega }\left( \int_{0}^{t}g(t-s)\left(|\nabla u(x,s)-\nabla
u(x,t)|\ +|\nabla u(x,t)|\right)\ ds\ \right) ^{2}\ dx.
\end{eqnarray}
Using the fact that 
\[
\int_{0}^{t}g(s)\ ds\leq 1-l 
\]
and Young's and H\"{o}lder's inequalities, we obtain, for any $\eta >0,$ 
\[
\int_{\Omega }\left( \int_{0}^{t}g(t-s)\left(|\nabla u(x,s)-\nabla u(x,t)|\
+|\nabla u(x,t)|\right)\ ds\ \right) ^{2}\ dx
\]
\begin{eqnarray*}
&=&\int_{\Omega }\left( \int_{0}^{t}g(t-s)\left(|\nabla u(s)-\nabla u(t)|\right)\right)
^{2}\ dx+\int_{\Omega }\left( \int_{0}^{t}g(t-s)|\nabla u(t)|\ ds\right)
^{2}\ dx \\
&&+2\int_{\Omega }\left( \int_{0}^{t}g(t-s)(|\nabla u(s)-\nabla u(t)|)\
ds\right) \quad \left( \int_{0}^{t}g(t-s)|\nabla u(t)|\ ds\right) \ dx
\end{eqnarray*}
\begin{eqnarray}
&\leq &(1+\eta )\int_{\Omega }\left( \int_{0}^{t}g(t-s)|\nabla u(t)|\
ds\right) ^{2}\ dx+\left(1+\frac{1}{\eta }\right)\int_{\Omega }\left(
\int_{0}^{t}g(t-s)|\nabla u(s)-\nabla u(t)|\ ds\right) ^{2}\ dx  \nonumber \\
&\leq &(1+\eta )(1-l)^{2}\Vert \nabla u(.,t)\Vert _{2}^{2}+\left(1+\frac{1}{\eta }%
\right)(1-l)(g\circ \nabla u)(t).
\end{eqnarray}
Substuting (3.6) in (3.5), we get 
\begin{eqnarray}
\int_{\Omega }\nabla u(x,t)\ .\int_{0}^{t}g(t-s)\nabla u(x,s)\ ds\ dx &\leq &%
\frac{1}{2}\left( 1+(1+\eta )(1-l)^{2}\right) \Vert \nabla u(.,t)\Vert
_{2}^{2}  \nonumber  \label{12} \\
&+&\frac{1}{2}\left(1+\frac{1}{\eta }\right)(1-l)(g\circ \nabla u)(t).
\end{eqnarray}
Combining (3.3), (3.4) and (3.7), we find 
\begin{eqnarray}
\Vert \nabla u(.,t)\Vert _{2}^{2} &\leq &c_{1}\delta E(t)-\frac{C_{\delta }}{%
c_{0}}E^{\prime }(t)  \nonumber  \label{13} \\
&+&\frac{1}{2}\left( 1+(1+\eta )(1-l)^{2}\right) \Vert \nabla u(.,t)\Vert
_{2}^{2}  \nonumber \\
&+&\frac{1}{2}\left(1+\frac{1}{\eta }\right)(1-l)(g\circ \nabla u)(t).
\end{eqnarray}
We then choose $0<\eta <$ $l(2-l)/(1-l)^{2},$ which makes $c_{2}=\frac{1}{2}%
\left( 1+(1+\eta )(1-l)^{2}\right) <1,$ and, therefore, (3.8) takes the form 
\[
\Vert \nabla u(.,t)\Vert _{2}^{2}\leq c_{1}\delta E(t)-\frac{C_{\delta }}{%
c_{0}}E^{\prime }(t)+c_{2}\Vert \nabla u(.,t)\Vert _{2}^{2}+c_{3}(g\circ
\nabla u)(t),
\]
where $c_{3}=\frac{1}{2}(1+\frac{1}{\eta })(1-l)$. This yields (3.2) with $%
c_{4}=\frac{c_{1}}{1-c_{2}}$ and $c_{5}=\frac{c_{3}}{1-c_{2}}$.\vspace{0.5cm}%
\newline
\textbf{Theorem 3.3}. \emph{Let }$u$\emph{\ be the solution of }$(1.1).$%
\emph{\ Then,} \emph{there exist strictly two positive constants }$\lambda _{0}$ and $%
\lambda _{1}$\emph{\ such that the energy satisfies, for all }$t\in \mathrm{I\hskip-2ptR}^{+}$\emph{, } 
\begin{equation}
E(t)\leq \lambda _{0} e^{-\lambda _{1}\int_{0}^{t}\xi (s)\ ds}\qquad \text{if }p=1,
\end{equation}
\begin{equation}
E(t)\leq \lambda _{0}\left( 1+\int_{0}^{t}\xi ^{2p-1}(s)\ ds\right) ^{\frac{-1}{2p-2}}\qquad \text{if }p>1.
\end{equation}
\emph{Moreover, if }$\xi $ \emph{and }$p$\emph{\ in }$(G_{2})$\emph{\ satisfy} 
\begin{equation}
\int_{0}^{+\infty }\left( 1+\int_{0}^{t}\xi ^{2p-1}(s)\ ds\right)
^{\frac{-1}{2p-2}}\ dt<+\infty ,
\end{equation}
\emph{then, for all }$t\in \mathrm{I\hskip-2ptR}^{+}$\emph{, }
\begin{equation}
E(t)\leq \lambda _{0} \left( 1+\int_{0}^{t}\xi ^{p}(s)\ ds\right) ^{\frac{-1}{%
p-1}} \qquad \text{if }p>1.
\end{equation}
\textbf{Remark 3.1}. Estimates (3.10) and (3.11) yield 
\begin{equation}
\int_{0}^{+\infty }E(t)\ dt<+\infty .
\end{equation}
\textbf{Proof.} From (3.1) and for any\ $\kappa >0,$ we have 
\begin{eqnarray*}
E^{\prime }(t) &\leq &0=-\kappa E(t)+\kappa E(t) \\
&\leq &-\kappa E(t)+\kappa \left( \frac{1}{2}(g\circ \nabla u)(t)+\frac{1}{2}%
\left(1-\int_{0}^{t}g(s)\ ds\right)\ \Vert \nabla u(.,t)\Vert _{2}^{2}\right)  \\
&\leq &-\kappa E(t)+\frac{\kappa }{2}(g\circ \nabla u)(t)+\frac{\kappa }{2}%
\Vert \nabla u(.,t)\Vert _{2}^{2}.
\end{eqnarray*}
Recalling Lemma 3.2, we get
\begin{eqnarray*}
E^{\prime }(t) &\leq &-\kappa E(t)+\frac{\kappa }{2}(g\circ \nabla u)(t) \\
&&+\frac{\kappa }{2}\left( c_{4}\delta E(t)-\frac{C_{\delta }}{c_{0}}%
E^{\prime }(t)+c_{5}(g\circ \nabla u)(t)\right)  \\
&\leq &-\kappa \left( 1-\frac{c_{4}}{2}\delta \right) E(t)-\frac{\kappa
C_{\delta }}{2c_{0}}E^{\prime }(t)+\frac{\kappa (1+c_{5})}{2}(g\circ \nabla
u)(t).
\end{eqnarray*}
Then we have 
\[
\left(1+\frac{\kappa C_{\delta }}{2c_{0}}\right)E^{\prime }(t)\leq -\kappa \left( 1-%
\frac{c_{4}}{2}\delta \right) E(t)+\frac{\kappa (1+c_{5})}{2}(g\circ \nabla
u)(t).
\]
By choosing $\delta $ small enough, we obtain, for two positive constants $\lambda $ and $\gamma $, 
\begin{equation}
E^{\prime }(t)\leq -\lambda E(t)+\gamma (g\circ \nabla u)(t).
\end{equation}
\underline{\textbf{Case of} $p=1$}. Multiplying (3.14) by $\xi (t)$ and
exploiting $(G_{2})$, we get 
\begin{eqnarray}
\xi (t)E^{\prime }(t) &\leq &-\lambda \xi (t)E(t)+\gamma (\xi g\circ \nabla
u)(t)  \nonumber  \label{26} \\
&\leq &-\lambda \xi (t)E(t)-\gamma (g^{\prime }\circ \nabla u)(t)  \nonumber
\\
&\leq &-\lambda \xi (t)E(t)-\gamma E^{\prime }(t).
\end{eqnarray}
We then set $L=(\xi +\gamma )E\sim E$ to obtain, from (3.15) and the fact
that $\xi ^{\prime }\leq 0$, 
\begin{equation}
L^{\prime }(t)\leq -\lambda \xi (t)E(t)\leq -\lambda _{1}\xi (t)L(t).
\end{equation}
A simple integration of (3.16) leads to 
\[
L(t)\leq Ce^{-\lambda _{1}\int_{0}^{t}\xi (s)\ ds}.
\]
This gives (3.9), by virtue of $L\sim E$.%
\vspace{0.5cm}\newline
\underline{\textbf{Case of} $p>1$}. To establish (3.10), we again consider
(3.14) and use Lemma 2.1 to get 
\[
\xi (t)E^{\prime }(t)\leq -\lambda \xi (t)E(t)+C\left( -E^{\prime
}(t)\right) ^{\frac{1}{2p-1}} .
\]
Multiplication of the last inequality by $\xi ^{\alpha }E^{\alpha }(t)$,
where $\alpha =2p-2>0$, gives 
\[
\frac{1}{\alpha +1}\xi ^{\alpha +1}\frac{d}{dt}E^{\alpha +1}(t)\leq -\lambda
\xi ^{\alpha +1}(t)E^{\alpha +1}(t)+c\left( \xi E\right) ^{\alpha }(t)\left(
-E^{\prime }(t)\right) ^{\frac{1}{\alpha +1}}.
\]
Use of Young's inequality, with $q=\alpha +1$ and $q^{*}=\frac{\alpha +1}{%
\alpha }$, yields, for any $\varepsilon >0,$ 
\begin{eqnarray*}
\frac{1}{\alpha +1}\xi ^{\alpha +1}\frac{d}{dt}E^{\alpha +1}(t) &\leq
&-\lambda \xi ^{\alpha +1}(t)E^{\alpha +1}(t)+C\left( \varepsilon \xi
^{\alpha +1}(t)E^{\alpha +1}(t)-C_{\varepsilon }E^{\prime }(t)\right)  \\
&=&-(\lambda -\varepsilon C)\xi ^{\alpha +1}(t)E^{\alpha
+1}(t)-C_{\varepsilon }E^{\prime }(t).
\end{eqnarray*}
We then choose $0<\varepsilon <\frac{\lambda }{C}$ and recall that $\xi
^{\prime }\leq 0$, to obtain, for $c_{6}>0,$ 
\[
\left( \xi ^{\alpha +1}E^{\alpha +1}(t)\right) ^{\prime }(t)\leq \xi
^{\alpha +1}\frac{d}{dt}E^{\alpha +1}(t)\leq -c_{6}\xi ^{\alpha
+1}(t)E^{\alpha +1}(t)-CE^{\prime }(t);
\]
which implies 
\[
\left( \xi ^{\alpha +1}E^{\alpha +1}+CE\right) ^{\prime }(t)\leq -c_{6}\xi
^{\alpha +1}(t)E^{\alpha +1}(t).
\]
Let $W=\xi ^{\alpha +1}E^{\alpha +1}+CE\sim E$. Then 
\[
W^{\prime }(t)\leq -C\xi ^{\alpha +1}(t)W^{\alpha +1}(t)=-C\xi
^{2p-1}(t)W^{2p-1}(t).
\]
Integrating over $(0,t)$ and using the fact that $W\sim E$, we obtain, for some $\lambda_0 >0$, 
\[
E(t)\leq \lambda_0\left(\int_{0}^{t}\xi ^{2p-1}(s)\ ds+1\right) ^{\frac{-1}{2p-2}};
\]
so (3.10) holds
\vspace{0.5cm}\newline
To establish (3.12), we put
\[
\eta (t)=\int_{0}^{t}\Vert \nabla u(t)-\nabla u(t-s)\Vert _{2}^{2} \ ds.
\]
Using Remark 3.1, we have
\begin{eqnarray*}
\eta (t) &\leq
&2\int_{0}^{t}\left( \Vert \nabla u(t)\Vert _{2}^{2}+\Vert \nabla
u(t-s)\Vert _{2}^{2}\right)\ ds \\
&\leq &\frac{4}{1-l}\int_{0}^{t}\left( E(t)+E(t-s)\right)\ ds \\
&=&\frac{8}{1-l}\int_{0}^{t}E(s)\ ds<\frac{8}{1-l}\int_{0}^{+\infty
}E(s)\ ds<+\infty .
\end{eqnarray*}
This implies that 
\begin{eqnarray}
\sup_{t\in \mathrm{I\hskip-2ptR}^{+}}\eta^{1-\frac{1}{p}} (t) <+\infty.
\end{eqnarray}
Assume that $\eta (t)>0$. Then, from (3.14), we find
\begin{eqnarray}
\xi (t)E^{\prime }(t) &\leq &-\lambda \xi (t)E(t)+\gamma \xi (t)(g\circ
\nabla u)(t) \\
&=&-\lambda \xi (t)E(t)+\gamma \frac{\eta (t)}{\eta (t)}\int_{0}^{t} \left(\xi
^{p}(s)g^{p}(s)\right)^{\frac{1}{p}}\Vert \nabla u(t)-\nabla u(t-s)\Vert
_{2}^{2}\ ds.  \nonumber
\end{eqnarray}
Applying Jensen's inequality (2.4) for the second term of the right-hand
side of (3.18), with 
\[
\Omega =[0,t],\quad f(s)=\xi ^{p}(s)g^{p}(s)\quad \mbox{and}\quad h(s)=\Vert \nabla u(t)-\nabla u(t-s)\Vert _{2}^{2} ,
\]
to get 
\[
\xi (t)E^{\prime }(t)\leq -\lambda \xi (t)E(t)+\gamma\eta (t)\left( \frac{1}{\eta
(t)}\int_{0}^{t}\xi ^{p}(s)g^{p}(s)\Vert \nabla u(t)-\nabla u(t-s)\Vert
_{2}^{2}\ ds\right) ^{\frac{1}{p}}.
\]
Therefore, using (3.17) we obtain 
\begin{eqnarray*}
\xi (t)E^{\prime }(t) &\leq &-\lambda \xi (t)E(t)+\gamma\eta^{1-\frac{1}{p}} (t)\left( \xi
^{p-1}(0)\int_{0}^{t}\xi (s)g^{p}(s)\Vert \nabla u(t)-\nabla u(t-s)\Vert
_{2}^{2}\ ds\right) ^{\frac{1}{p}} \\
&\leq &-\lambda \xi (t)E(t)+C(-g^{\prime }\circ \nabla u)^{\frac{1}{p}}(t),
\end{eqnarray*}
and then
\begin{eqnarray}
\xi (t)E^{\prime }(t) \leq -\lambda \xi (t)E(t)+C(-E^{\prime }(t))^{\frac{1}{p}} .
\end{eqnarray}
If $\eta (t)=0$, then $s \rightarrow\nabla u (s)$ is a constant function on $[0,t]$. Therefore 
\[
(g\circ \nabla u)(t)=0,
\]
and hence we have, from (3.14), 
\[
E^{\prime }(t) \leq -\lambda E(t),
\]
which implies (3.19). 
\vspace{0.5cm}\newline
Now, multiplying (3.19) by $\xi ^{\alpha }(t)E^{\alpha }(t)$, for $\alpha =p-1$, and
repeating the same computations as in above, we arrive at, for some $\lambda_0 >0$, 
\[
E(t)\leq \lambda_0 \left( \int_{0}^{t}\xi ^{p}(s)\ ds+1\right) ^{\frac{-1}{p-1}}.
\]
This completes the proof of our main result.\vspace{0.5cm}

The following examples illustrate our result and show the optimal decay rate
in the polynomial case:\vspace{0.5cm}\newline
\textbf{Example 3.1.} Let $g(t)=a(1+t)^{-\nu },$ where $\nu >2,$ and $a>0$ so that 
\begin{equation}
\int_{0}^{+\infty }g(t)\ dt<1. 
\end{equation}
We have 
\[
g^{\prime }(t)=-a\nu (1+t)^{-\nu -1}=-b\left( a(1+t)^{-\nu}%
\right) ^{\frac{\nu +1}{\nu }} ,
\]
where $b=\nu a^{-\frac{1}{\nu }}$. Then $(G_2)$ holds with $\xi (t)=b$ and $p=\frac{\nu +1}{\nu }\in (1,\frac{3}{2})$. Therefore (3.11) yields 
\[
\int_{0}^{+\infty }\left( b^{2p-1}t+1\right) ^{\frac{-1}{2p-2}}\ dt<+\infty , 
\]
and hence, by (3.12), we get 
\[
E(t)\leq C(1+t)^{\frac{-1}{p-1}}=C(1+t)^{-\nu },
\]
which is the optimal decay.\vspace{0.5cm}\newline
\textbf{Example 3.2}. Let $g(t)=ae^{-(1+t)^{\nu }}$, where $0<\nu \leq 1$, and $a>0$ is chosen so that (3.20) holds. Then 
\[
g^{\prime }(t)=-a\nu (1+t)^{\nu -1}e^{-(1+t)^{\nu }}.
\]
Therefore $(G_2)$ holds with $p=1$ and $\xi (t)=\nu (1+t)^{\nu -1}$. Consequently, we can use (3.9) to deduce 
\[
E(t)\leq Ce^{-\lambda (1+t)^{\nu }}.
\]

\section*{Acknowledgment}

The authors thank KFUPM and Lorraine-Metz university for their continuous
support.  This work has been finalized during the visit of the third author
to KFUPM in December 2016 and during the scholarship of the first author in
Lorraine-Metz. This work has
been partially funded by KFUPM under Project \# IP152-Math 212.

\section{References}

\begin{enumerate}
\item  Berrimi S. and Messaoudi S.A., A decay result for a quasilinear
parabolic system, Progress in Nonlinear Differential Equations and their
Applications \textbf{53} (2005), 43-50.

\item  Da Prato G. and Iannelli M., Existence and regularity for a class of
integro-differential equations of parabolic type, J. Math. Anal. Appl.%
\textit{\ }\textbf{112} (1985), 36-55.

\item  Englern H., Kawohl B. and Luckhaus S., Gradient estimates for
solutions of parabolic equations and systems, J. Math. Anal. Appl.\textit{\ }%
\textbf{147} (1990), 309-329.

\item  Liu G. and Chen H., \textit{Global and blow-up of solutions for a
quasilinear parabolic system with viscoelastic and source terms}, Math.
Methods Appl. Sci. \textbf{37 }(2014), 148-156.

\item  Messaoudi S.A. and Tellab B., A general decay result in a quasilinear
parabolic system with viscoelastic term, Applied Mathematics Letters \textbf{%
25} (2012), 443-447.

\item  Messaoudi S.A. and Al-Khulaifi W., General and optimal decay for a
quasilinear viscoelastic equation, Applied Mathematics Letters \textbf{66}
(2017), 16-22.

\item  Nakao M. and Ohara Y., Gradient estimates for a quasilinear parabolic
equation of the mean curvature type, \textit{J. Math. Soc. Japan}\textbf{\
48 \# 3} (1996), 455-466.

\item  Nakao M. and Chen C., Global existence and gradient estimates for the
quasilinear parabolic equations of $m$-Laplacian type with a nonlinear
convection term,\textit{\ }Journal of Differential Equations \textbf{162}
(2000), 224-250\textit{.}

\item  Nohel J.A., Nonlinear Volterra equations for the heat flow in
materials with memory, Integral and functional differential equations,
Lecture notes in Pure and Applied Mathematics, Marcel Dekker Inc. 1981%
\textit{.}

\item  Pucci P. and Serrin J., Asymptotic stability for nonlinear parabolic
systems, Energy methods in continuum mechanics, Kluwer Acad. Publ.,
Dordrecht, 1996.

\item  Yin H.M., On parabolic Volterra equations in several space
dimensions, \textit{SIAM J. of Mathematical Analysis. }\textbf{22} (1991),
1723-1737.
\end{enumerate}

{}

\end{document}